\documentclass[11pt,a4paper,reqno]{amsart}
\usepackage{amsthm, amssymb, amscd} 
\usepackage{enumerate,array}
\usepackage{oppi}
\usepackage{extpfeil}
\newtheorem{theorem}{Theorem}[section]
\newtheorem{lemma}[theorem]{Lemma}

\newtheorem{proposition}[theorem]{Proposition}
\theoremstyle{definition}

\newtheorem{question}[theorem]{Question}
\newtheorem{conjecture}[theorem]{Conjecture}

\theoremstyle{remark}
\newtheorem{remark}[theorem]{Remark}

\numberwithin{equation}{section}
\newcommand\dont{{\widetilde \dD_n^0}}
\begin{document} 
%

\title[Spectral Gap of $W(B_n)$]{%
On the spectral gap of some Cayley graphs on the Weyl group $W(B_n)$
}
%
%
%
%
\author[Filippo Cesi]{Filippo Cesi}
\address{%
    Filippo Cesi\hfill\break
    \indent Dipartimento di Fisica\hfill\break
    \indent Universit\`a di Roma ``La Sapienza", Italy.
}
\email{filippo.cesi@roma1.infn.it}



%
%
%
\begin{abstract}
The Laplacian of a (weighted) Cayley graph on the Weyl
group $W(B_n)$ is a 
$N\times N$ matrix
with $N = 2^n n!$ equal to the order of the group.
We show that for a class of (weighted) generating sets,
its spectral gap (lowest nontrivial eigenvalue),
is actually equal to 
the spectral gap of a
$2n \times 2n$ matrix associated to a $2n$-dimensional
permutation representation of $W_n$.
This result can be viewed as an extension to $W(B_n)$
of an analogous result valid for the symmetric group,
known as ``Aldous' spectral gap conjecture'', proven in
2010 by Caputo, Liggett and Richthammer.
\end{abstract}
%
%
\maketitle
\thispagestyle{empty}
%
%
%
%
\vspace{5mm}
\section{Introduction} 
\label{sec:intro}

\noindent
Let $G$ be a finite group with complex group algebra $\bC G$.
If $w = \sum_{g\in G} w_g \, g$ is an element of the group
algebra such that all coefficients $w_g$ are real, nonnegative and symmetric,
\ie\ $w_{g^{-1}} = w_g$, 
we denote with $\Cay(G,w)$ 
the weighted Cayley graph
whose vertices are the elements of $G$
and whose (undirected) edges
are the pairs $\{g,hg\}$ with $g,h \in G$. Each edge $\{g,hg\}$
carries a \textit{weight} equal to $w_h$.
The Laplacian of $\Cay(G,w)$
is a linear operator acting
on functions $f:G\to \bC$ as
\begin{align*}
  & \left[ \D_{\Cay(G,w)} f\right](g) = \sum_{h\in G} w_h \bigl( f(g) - f(hg) \bigr)
  \,.
\end{align*}
Since (weighted) Cayley graphs are regular, the Laplacian is
strictly related to the (weighted) adjacency matrix $A_{\Cay(G,w)}$,
namely
\begin{align*}
  \D_{\Cay(G,w)} = \Bigl( \sum_{g\in G} w_g \Bigr) \id_N - A_{\Cay(G,w)}
  \,,
\end{align*}
where $N$ is the order of $G$ and $\id_N$ is the $N\times N$
identity matrix.
The lowest eigenvalue of the Laplacian is trivially zero
with constant eigenvector.
The \textit{spectral gap} of $\Cay(G,w)$ is defined as the second
lowest eigenvalue of the associated Laplacian and it is denoted with
$\psi_G(w)$. It is strictly positive if and only if
the support of $w$ generates $G$, that is if $\Cay(G,w)$ is connected.

\medno
Taking into account the symmetry of $w$, we can rewrite the Laplacian as
\begin{align}
  \label{eq:L}
  \D_{\Cay(G,w)} = \sum_{h\in G} w_h \left[\id_N - \dL(h) \right]
  \,,
\end{align}
where $\dL$ is the left regular representation of $G$ acting
on functions on $G$ as
\begin{align}
  \label{eq:left}
  [\dL(h) f](g) = f(h^{-1} g)
  \,.
\end{align}
Formula \eqref{eq:L} suggests the following generalization:
Given a representation $\dR$ of $G$ 
on the $d$-dimensional complex vector space $V$, 
and given $w \in \bC G$, following \cite{cesiOctopus},
we define the \textit{representation Laplacian} $\D_G(w,\dR)$
as the linear operator on $V$ given by
\begin{align}
  \label{eq:lapl}
  &\D_G(w,\dR) := \sum_{h\in G} w_h \, \left[\id_V - \dR(h)\right] 
  && w_h \in \bC
  ,
\end{align}
where $\id_V$ is the identity on $V$.
To the pair $(w,\dR)$ we also associate a
\textit{spectral gap}, denoted by $\psi_G(w,\dR)$, which is again
the smallest
nontrivial eigenvalue of $\D_G(w,\dR)$ 
(see Section \ref{sec:notation} for a precise
definition).
Thus the Laplacian of the Cayley graph is a special
case of the representation Laplacian, and we can write
\begin{align}
  \label{eq:psiL}
  &\D_{\Cay(G,w)} = \D_G(w,\dL)
  &&\psi_G(w) = \psi_G(w,\dL) 
  \,.
\end{align}
In this paper we pursue the general idea of \cite{cesiOctopus}
that, although $\dL$ contains all irreducible representations
of $G$, in some interesting cases
it is possible to pinpoint those representations
which are ``responsible'' for the spectral gap of the Cayley graph.
These representations can then replace $\dL$ in \eqref{eq:psiL},
with the advantage of having to deal with a possibly
\textit{much smaller} matrix
than $\D_{\Cay(G,w)}$.

The most important result so far in this direction is the proof 
of the so called \textit{Aldous' spectral gap conjecture},
concerning the symmetric group.
After several partial successes
in a series of papers \cite{DiSh}, \cite{FOW},  \cite{Bac}, \cite{HaJu},
\cite{KoNa}, \cite{Mor}, \cite{StCo}, \cite{Ces1}, \cite{Diek} 
spanning about 25 years, a general proof was finally given in
\cite{CaLiRi}.
While the original formulation was given in a probabilistic
framework \cite{Ald}, the statement can be translated
as follows (see \cite{Ces1} for more details on the equivalence):
let $\fS_n$ be the symmetric group on $\{1,\ldots,n\}$, and
let $T_n$ be the set of all transpositions in $\fS_n$.

\begin{theorem}\label{thm:ald}
(Aldous' spectral gap conjecture, proven in \cite{CaLiRi}).
Let $w\in \bC \fS_n$ be given by
\begin{equation*}
  w = \sum_{(ij)\in T_n} b_{ij}\, (ij)
  \,,
\end{equation*}
where $(ij)$ is the
transposition which exchanges $i$ and $j$, and $b_{ij} \ge 0$.
Then
\begin{equation}
  \label{eq:D}
  \psi_G(w) = \psi_G(w, \dD^0_n) \,,
\end{equation}
where $\dD^0_n$ is the $n$-dimensional \textit{defining} representation
of $\fS_n$, associated with the natural action of
$\fS_n$ on the set $\{1,2,\ldots,n\}$.
\end{theorem}

\smallno
Thanks to this theorem, if $w$ is supported on transpositions,
in order to find
the spectral gap of the Laplacian of $\Cay(\fS_n,w)$, which
is a $n!\times n!$ matrix, all one has to do
is to find the smallest nontrivial eigenvalue of a $n\times n$
matrix.

The spectral gap of Cayley graphs on the symmetric or alternating groups
has been computed also in some particular cases
where the generators are not transpositions.
Examples are the initial reversal graph \cite{Ces3},
the (complete, extended) alternating group graph \cite{huang2017adjacency} and
the substring reversal graph \cite{chung2017spectral}.
In \cite{parzanchevski2018aldous} the authors prove that
if $w = \sum_{\pi\in [\a]} \pi$, where $[\a]$ is a conjugacy class
of $\fS_n$, then, for $n$ large enough, the spectral gap of $w$ is 
associated with one of eight low-dimensional representations.
A weaker statement is also proven for the more general case
where $w$ is invariant under conjugation.

In \cite{cesiOctopus} we described a general strategy
for proving results analogous to Theorem \ref{thm:ald}
for arbitrary finite groups
and we gave a slightly simpler proof
based on this point of view.
In particular, in this strategy, it appears that
the representations ``responsible'' for the spectral gap
are the nontrivial irreducible representations of $\fS_n$ which,
when restricted to $\fS_{n-1}$, contain the trivial
representation.
In the case of the symmetric group there is a unique such
representation, namely the one associated with the partition
$(n-1,1)$. This representation, apart from a trivial summand,
is equivalent to $\dD^0_n$ (see \eqref{eq:d0}).

In this paper we apply the idea of \cite{cesiOctopus} 
to the Weyl/Coxeter group
$W(B_n)$ associated with the $B_n$ (or $C_n$) root system,
also called the \textit{hyperoctahedral group}.
For simplicity we let $W_n = W(B_n)$.
There are several equivalent ways to define this group.
One possible realization%
is as the subgroup of $GL(n,\bC)$ consisting 
of all $n\times n$ matrices which
have exactly one non-zero entry in each row and each column,
and this non-zero entry is either $1$ or $-1$.
We have thus a natural embedding $\fS_n \inje{} W_n$
where the symmetric group is the subgroup of all matrices with 
nonnegative entries.
The group $W_n$ can also be described as a group of
\textit{signed permutations}. This leads to another embedding
\begin{align*}
  \dP_n : W_n \inje{} \fS_{2n}
  \,,
\end{align*}
where $\dP_n$ is a $2n$-dimensional faithful permutation
representation described in more details in Section \ref{sec:wn}.

\medno
We can now state the main result of this paper:
for $A\sset \{1,\ldots,n\}$,
let $s_{\{i\}} \in W_n$ be the diagonal matrix%
\footnote{In this introduction we are implicitly using the
\textit{defining representation} for describing the
elements of $W_n$. See \eqref{eq:Dsi}.
This particular realization of $W_n$ will be denoted with
$\widetilde W_n$ in Section \ref{sec:wn}.}
\begin{equation}
  s_{\{i\}} = \diag(1,\ldots,1,-1,1,\ldots,1) 
  \,,
\end{equation}
where the unique $-1$ occurs in the $i^{th}$ place.  
Then we have:

\begin{theorem}\label{thm:main}
Let $w\in \bC W_n$ be given by
\begin{align}
  \label{eq:wmain}  
  &w = \sum_{i=1}^n a_i\, s_{\{i\}} + \sum_{(ij)\in T_n} b_{ij}\, (ij)
  &&a_i\ge0,\; b_{ij} \ge 0
  \,.
\end{align}
Then 
\begin{align}
  \label{eq:main}
  \psi_{W_n}(w) = \psi_{W_n}(w, \dP_n )
  \,.
\end{align}
\end{theorem}

\medno
In our approach the representation $\dP_n$ appears in
\eqref{eq:main} for the same reason that the defining
representation of the symmetric group appears in \eqref{eq:D},
that is $\dP_n$ contains all irreducible representations
of $W_n$ which, when restricted to $W_{n-1}$, contain
the trivial representation. 

There is another result which is worth mentioning
since it has an interesting overlap with Theorem \ref{thm:main}.
In \cite{Kas} it is proved that if $(G,S)$ is a finite
Coxeter system, then both the spectral gap and the Kazhdan
constant are determined by
the defining
representation.
In our notation, this implies that if $G=W(B_n)$ and
if $w$ is the sum of a Coxeter generating set, that is
\begin{equation}
  \label{eq:cox}
  w = s_{\{1\}} + \sum_{i=1}^{n-1} (i, i+1)
\end{equation}
the $\psi_{W_n}(w) = \psi_{W_n}(w, \dD_n)$.
This result covers every finite Coxeter group,
but applies (essentially) to
only one element of the group algebra, namely $w = \sum_{s\in S} s$,
where $S$ is a Coxeter generator for $G$.
The approach used in the proof of this theorem is very different
from ours, and in Section 5 of \cite{cesiOctopus} we
explain why it is unlikely that it could be effective
for dealing with more general elements $w$ of the group algebra.
Since the defining representation $\dD_n$ of $W_n$
is a subrepresentation of $\dP_n$, Kassabov's result
is stronger than ours for $w$ of the form \eqref{eq:cox}.
Nevertheless, \eqref{eq:main} is optimal for a
general $w$ of type \eqref{eq:wmain}.
In Section \ref{sec:rem} we show, in fact, that the theorem 
is (in general) false if we try to improve it by 
replacing $\dP_n$ with a subrepresentation.
We also discuss a possible generalization of Theorem \ref{thm:main}.

\medno

\section{The representation Laplacian and its spectral gap} 
\label{sec:notation}

\noindent
If $G$ is a finite group, $\rep(G)$ denotes the set
of all finite-dimensional complex representations of $G$,
while $\Irr(G)$ is the set of all equivalence classes of
irreducible representations.
By Maschke's theorem, 
we have for each $\dR \in \rep(G)$,
\begin{equation}\label{eq:Y}
  \dR \cong \bigoplus_{\dT \in \Irr(G)} \mu(\dR, \dT) \, \dT \,,
\end{equation}
where $\mu(\dR, \dT)$ is a nonnegative integer
called the \textit{multiplicity} of $\dT$ in $\dR$.
If $\dR$ is a representation of $G$ on the complex vector space $V$,
$V^{G,\dR}$ stands for
the subspace of all invariant vectors
\begin{equation*}
  V^{G,\dR} = \{ v \in V : \dR(g) v = v,\ \forall g\in G \}.
\end{equation*}
By definition we have
\begin{align}
  \label{eq:mult}
  \dim V^{G,\dR} = \mu(\dR,\dI)
  \,,
\end{align}
where $\dI$ is the
one-dimensional trivial representation.
An eigenvalue $\l$ of the representation Laplacian 
$\D_G(w,\dR)$, defined in \eqref{eq:lapl}, will be called \textit{trivial}
if its corresponding eigenspace consists entirely
of invariant vectors $v \in V^{G,\dR}$. 

If $w \in \bC G$, the \textit{support} of $w$ is defined as
\begin{align*}
  \supp w = \{ g \in G : w_g \ne 0 \}
  \,.
\end{align*}
We introduce a canonical involution in the group algebra $\bC G$ as
\begin{equation*}
  w = \sum_{g\in G} w_g \, g \too w^* := 
  \sum_{g\in G} \ol w_g \, g^{-1} 
  \,.
\end{equation*}
An element $w\in \bC G$ is called \textit{symmetric} if $w = w^*$, 
and it is called \textit{positive} if $w_g \ge 0$ for all $g\in G$.
We let
\begin{align*}
  \bC G^{(s)} &= \{ w \in \bC G : \text{$w$ is symmetric} \} 
  \\
  \bR_+ G^{(s)} &= \{ w\in \bC G : \text{$w$ is symmetric and positive} \}.
\end{align*}
It is easy to prove \cite[Sect. 2]{cesiOctopus} that
$\D_G(w,\dR)$ has real eigenvalues if $w$ is symmetric and
real nonnegative eigenvalues if $w$ is symmetric and positive.
If $w$ is symmetric we can label the eigenvalues of $\D_G(w,\dR)$,
with possible repetitions according to their multiplicity,
in nondecreasing order as
\begin{align*}
  \l_1( \D_G(w,\dR) ) \le \l_2(\D_G(w,\dR)) \le \cdots 
  \le \l_s( \D_G(w, \dR) )
  \,,
\end{align*}
where $s$ id the degree (or dimension) of the representation $\dR$.
If $w \in \bR_+ G^{(s)}$,
we define the \textit{spectral gap}
of the pair $(w, \dR)$ as
\begin{align}
  \label{eq:sgdef}
  \psi_G(w, \dR ) &:= \min \{ \l \in \spec \D_G(w,\dR) :
  \text{$\l$ is nontrivial} \}
  ,
\end{align}
with the convention that $\min \emp = +\oo$.
If $t = \dim V^{G,\dR}$, then $\D_G(w,\dR)$ has exactly $t$
trivial eigenvalues, thus, thanks to \eqref{eq:mult}, we have%
\footnote{unless $\dR$ is a multiple of $\dI$ in which case
all eigenvalues are trivial, and thus the spectral gap is equal to $+\oo$.}
\begin{align}
  \label{eq:gapt}
  &\psi_G(w,\dR) = \l_{t+1}\left( \D_G(w,\dR) \right)
  &&\text{where $t = \mu(\dR,\dI)$}.
\end{align}
The \textit{spectral gap of $w$} is defined by minimizing
over representations
\begin{align}
  \label{eq:minrep1}
  \psi_G(w) &= \inf \bigl\{ \psi_G(w,\dR) : 
  \dR \in \rep(G)
  \bigr\} \,.
\end{align}
From \eqref{eq:Y} it follows that 
\begin{align}
  \label{eq:Y1}
  \spec \D_G( w,\dR) &= \bigcup_{\dT\in \Irr(G) :\, \mu(\dR,\dT)>0} 
  \spec \D_G( w,\dT ) 
\end{align}  
which implies
\begin{align}
  \label{eq:Y2}
  \psi_G(w,\dR) &= \min\{ \psi_G(w,\dT) : \dT\in \Irr(G),\ \mu(\dR,\dT) > 0 \}
  \,.
\end{align}
By consequence, in \eqref{eq:minrep1} we can just
consider irreducible representations, so
\begin{align}
  \label{eq:minrep2}
  \psi_G(w) &= \min \bigl\{ \psi_G(w,\dR) : 
  \dR \in \Irr(G)
  \bigr\} \,.
\end{align}
Let $\dL$ be the left regular representation of $G$
defined in \eqref{eq:left}. Since
\begin{equation}
  \label{eq:LT}
  \dL = \bigoplus_{\dT\in \Irr(G)} f_\dT\, \dT
  \,,
\end{equation}
where $f_\dT$ is the degree od $\dT$, we have $\mu(\dL,\dI)=1$.
Therefore 
$\psi_G(w) = \psi_G(w,\dL) = \l_2( \D_G(w, \dL))$.
This shows that definitions 
\eqref{eq:minrep1} and \eqref{eq:minrep2} for the quantity
$\psi_G(w)$ actually agree with the definition given in
Section \ref{sec:intro}
as the second lowest eigenvalue of the Cayley graph
$\Cay(G,w)$.

\section{Groups $\fS_n$, $W(B_n)$ and their representations}
\label{sec:wn}

\noindent
In this section we review
some more or less well known facts about the symmetric group,
the Weyl (or Coxeter) group $W_n:=W(B_n)$ and their representations.
For more details we refer the reader to \cite{GecPfe}, 
\cite{JaKe}, \cite{GeiKin}.

\noindent
A \textit{partition} of $n$
is a nonincreasing sequence $\a = (\a_1, \a_2, \ldots, \a_r)$
of positive integers  such that 
$\sum_{i=1}^r \a_i = n$. The \textit{size} of $\a$ is defined
as $|\a|=\sum_{i=1}^r \a_i$.
We write $\a\partit n$ if $\a$ is a partition of $n$.
The irreducible representations of $\fS_n$ are indexed (modulo equivalence) by
the partitions of $n$. If $\a\partit n$, we denote with $[\a]$
the corresponding irreducible representation of $\fS_n$.

The group $W_n$ can be realized as the set of all
pairs $(\h, \pi)$ with $\h \in \{0,1\}^n$ and $\pi \in \fS_n$
with product
\begin{align*}
  (\h,\pi) \cdot (\z, \s) = ( \h + \z\circ \pi^{-1}, \pi \s)
\end{align*}
where the sum of two elements of $\{0,1\}^n$ is componentwise
$\bmod\ 2$, and elements of $\{0,1\}^n$ are identified with
functions from $\{1,\ldots,n\}$ to $\{0,1\}$.
Observe that 
$(\h,\pi) = (\h, \idu_{\fS_n}) \cdot (0,\pi) = (0,\pi) \cdot (\h\circ \pi, \idu_{\fS_n})$.
Consider the two subgroups
\begin{align*}
  N_n &:= \{ (\h, \idu_{\fS_n}) : \h \in \{0,1\}^n \} \iso (\bZ/2\bZ)^n 
  \\
  S_n &:= \{ ( 0, \pi ) : \pi \in \fS_n \} \iso \fS_n
  \,.
\end{align*}
$N_n$ is a normal subgroup of $W_n$
and $W_n$ can be written as a semidirect product
\begin{align*}
  W_n = N_n \rtimes S_n \iso (\bZ/2\bZ)^n \rtimes \fS_n
  \,. 
\end{align*}
The irreducible representations of $W_n$ are indexed by 
ordered pairs of partitions $(\a,\b)$ such that 
$|\a|+|\b|=n$. We denote with $[\a,\b]$ the irreducible
representation corresponding to $(\a,\b)$.
We denote with $\dT^\a$ and $\dT^{(\a,\b)}$ some specific
(but arbitrary) choice of representations
in the equivalence classes $[\a]$ and $[\a,\b]$ respectively.
Given a pair of partitions $(\a,\b)$ with $|\a|=k$ and $|\b|=n-k$,
the representation $[\a,\b]$ can be obtained \cite[Sect. 2]{GeiKin} as an induced
representation as%
\footnote{we use the same notation as \cite{GecPfe}, while
in \cite{GeiKin} the order of $[\a,\b]$ is reversed.}
\begin{align}
  \label{eq:abind}
  [\a,\b] \cong \bigl( \dU_k \otimes [\a] \otimes [\b] \bigr)
  \indR_{N_n \times S_k \times S_{n-k}}^{W_n}
  \,,
\end{align}
where $\dU_k$ is the one-dimensional representation
of $N_n$ given by
\begin{align*}
  &\dU_k( \h ) = (-1)^{\card\{ i \in \{k+1,\ldots,n\} \,:\, \h_i = 1 \} }
  &&\h \in \{0,1\}^n
  ,\ k\in\{0,\ldots,n\}
  \,.
\end{align*}
In particular, when $k=n$ (and thus $\b=\emp$), we have that 
$\dU_n$ is the trivial representation and
$[\a,\emp]$
is the pullback of the $[\a]$ representation of $S_n \cong W_n/N_n$,
that is%
\begin{align}
  \label{eq:emp}
  &\dT^{(\a,\emp)}( \h, \pi ) = \dT^\a(\pi)
  &&(\h,\pi)\in W_n
  \,.
\end{align}
The trivial representation of $W_n$ is given by $\dI_n = [ (n), \emp ]$.

\medno
\textit{Branching rules}.
An irreducible representation of a finite group is in general
no longer irreducible when restricted to a subgroup,
but it can be expressed as a direct sum of irreducible representations
of the subgroup. 
The \textit{branching rule} $\fS_n \to \fS_{n-1}$
is 
\cite[Sect. 6.1.8]{GecPfe}
\begin{align}
  \label{eq:br}
  [\a] \resR^{\fS_n}_{\fS_{n-1}} = \bigoplus_{\b \in \a^-} \, [\b]
  \qquad \a\partit n
\end{align}
where, if $\a=(\a_1, \ldots,\a_r)$, 
$\a^-$ is defined as the collection of all sequences
of the form
\begin{equation*}
  (\a_1, \ldots, \a_{i-1}, \a_i -1, \a_{i+1}, \ldots, \a_r)
\end{equation*}
\textit{which are partitions of} $n-1$. 
The \textit{branching rule} $W_n \to W_{n-1}$ is
\cite[Sect. 6.1.9]{GecPfe}
\begin{align}
  \label{eq:bra}
  [ \a, \b ] \resR^{\cW_n}_{\cW_{n-1}} =
  \bigoplus_{\g \in \a^-} [\g, \b] 
  \oplus
  \bigoplus_{\g \in \b^-} [\a, \g] 
  \,.
\end{align}

\medno
\textit{The defining representation of $\fS_n$}.
Let $\dD_n^0$ be the defining $n$-dimensional representation
of $\fS_n$ with matrix elements
\begin{align}
  \label{eq:D0nij}
  &[\dD_n^0(\pi)]_{ij} = \d_{i,\pi(j)}
  &&\pi \in \fS_n
  \,.
\end{align}
This representation is not irreducible, but it can be decomposed
as
\begin{align}
  \label{eq:d0}
  \dD_n^0 = (n) \oplus (n-1,1)
  \,.
\end{align}

\medno
\textit{The defining representation of $W_n$}.
We let $\dD_n$ be the $n$-dimensional defining representation
of $W_n$ given by
\begin{align}
  \label{eq:Dnij}
  \left[\dD_n(\h,\pi)\right]_{ij} = (-1)^{\h_i} \, \d_{i, \pi(j)}
  \,.
\end{align}
This is a faithful representation, hence $W_n$ is
isomorphic to the image of $\dD_n$
which is the group $\widetilde W_n$
of all $n\times n$ matrices which
have exactly one non-zero entry in each row and each column,
and this non-zero entry is either $1$ or $-1$.
The normal subgroup $N_n$ is mapped to the subgroup
of the diagonal matrices of $\widetilde W_n$,
while the restriction of $\dD_n$ to $S_n$ is just
the $n$-dimensional defining representation of $\fS_n$.
It follows from \eqref{eq:abind} (see also \cite[Proposition 5.5.7]{GecPfe}
for a more general statement)
that $\dD_n$ is irreducible and that, in particular,
\begin{align}
  \label{eq:Dn}
  \dD_n \iso [ (n-1) , (1) ]
  \,.
\end{align}

\medno
\textit{The representation $\widetilde \dD_n^0$.}
Since $N_n$ is normal in $W_n$, every representation $\dR$
of the quotient $W_n/N_n \iso \fS_n$ can be pulled back (or lifted)
to a representation $\widetilde \dR$ of $W_n$ letting
\begin{align}
  \label{eq:lift}
  &\widetilde\dR( \h, \pi ) = \dR( \pi )
  &&(\h,\pi) \in W_n
  \,.
\end{align}
Furthermore $\widetilde \dR$ is irreducible
if and only if $\dR$ is.
We define $\widetilde \dD_n^0$ as the pullback of the
defining $n$-dimensional representation of $\fS_n$. Its matrix elements
are then
\begin{align}
  \label{eq:D0nija}
  [\widetilde \dD_n^0( \h, \pi )]_{ij} = \d_{i, \pi(j)}
  \,.
\end{align}
From \eqref{eq:d0} and \eqref{eq:emp} it follows that
\begin{align}
  \label{eq:D0n}
  \widetilde \dD_n^0 = [ (n), \emp ] \oplus [ (n-1,1), \emp ]
  \,.
\end{align}

\medno
\textit{The permutation representation $\dP_n$}.
Let $X_n =\{-n, \ldots, -1 \} \cup \{1, \ldots, n\}$
and consider the (left) group action of $W_n$ on $X_n$
given by
\begin{align}
  \label{eq:act}
  &(\h,\pi) k = (-1)^{ (\h\circ \pi)(|k|)} \, \sgn(k) \,\pi(|k|)
  &&k\in X_n
  \,.
\end{align}
We define $\dP_n$ as the $2n$-dimensional
permutation representation associated with this action. 
$\dP_n$ acts on the complex vector space
\begin{align*}
  V_n = \bC X_n := \{ (x_i)_{i\in X_n} : x_i \in \bC \}
  \,.
\end{align*}
If $(e_i)_{i\in X_n}$ is the canonical basis of $V_n$, such that
\begin{align*}
  \sum_{i\in X_n} x_i e_i = (x_{-n}, \ldots, x_{-1}, x_1, \ldots, x_n )
  \,,
\end{align*}
the representation matrices are determined by the equalities
\begin{align*}
  &\dP_n( g ) e_k = e_{ g k } 
  && g \in W_n,\ k\in X_n
  \,,
\end{align*}
where $gk$ is given by \eqref{eq:act}.
Therefore the matrix elements of $\dP_n$ are given by
\begin{align*}
  [\dP_n(\h,\pi)]_{ij} =
  \begin{cases}
  1 & \text{if $\pi(|j|) = |i|$ and $\sgn(j) = \sgn(i) \, (-1)^{\h_{|i|}}$} \\
  0 & \text{otherwise.}
  \end{cases}
\end{align*}
This representation is also faithful, so $W_n$ is isomorphic to
the image of $\dP_n$ which consists of the set of all permutations $\pi$ of $X_n$
such that $\pi(-k) = -\pi(k)$ for each $k\in X_n$ (the so called
\textit{signed permutations}).

In the following proposition we find the irreducible components
of $\dP_n$.

\begin{proposition}\label{thm:Pn}
We have
\begin{align}
  \label{eq:p+-}
  \dP_n = \dD_n \oplus \widetilde \dD_n^0 =
  [ (n-1), (1) ]
  \oplus 
  \dI_n \oplus [ (n-1,1) , \emp ] 
  \,.
\end{align}
\end{proposition}

\Pro\
Let $V_n = \bC X_n$ and $(e_i)_{i\in X_n}$ be as above, and let
\begin{align*}
   &e^+_k = e_k + e_{-k}
  &&e^-_k = e_k - e_{-k}
  && k = 1, \ldots, n
  \,.
\end{align*}
Let $V^+_n$ ($V^-_n$) be the subspace of $V_n$ spanned by
$(e^+_k)_{k=1}^n$ ($(e^-_k)_{k=1}^n$).
In other words $V_n^+$ is the subspace of the ``even''
vectors such that $x_{-i} = x_i$, while $V^-_n$ is the
subspace of the odd vectors. 

Let $g = (\h,\pi) \in W_n$.
The action defined in \eqref{eq:act} satisfies $g(-k) = - g(k)$.
By consequence we have, for $k=1,\ldots,n$,
\begin{align}
  \label{eq:Pe+}
  \dP_n(g) e^+_k 
  = e_{gk} + e_{g(-k)}
  = e_{gk} + e_{-gk}
  = e^+_{|gk|}
  = e^+_{\pi(k)}
\end{align}
and
\begin{equation}
\begin{split}
  \label{eq:Pe-}
  \dP_n(g) e^-_k 
  &= e_{gk} - e_{g(-k)}
  = e_{gk} - e_{-gk}
  = \sgn(gk) \, e^-_{|gk|}
  \\
  &=
  (-1)^{(\h\circ\pi)(k)} \, e^-_{\pi(k)} 
  \,.
\end{split}
\end{equation}
It follows
that 
both $V^+_n$ and $V^-_n$ are invariant under $\dP_n(g)$,
thus we have a direct sum decomposition
\begin{align*}
  &\dP_n = \dP^+_n \oplus \dP^-_n
  &&V_n = V^+_n \oplus V^-_n
  &&\dP_n^\pm := \dP_n \restriction_{V^\pm_n}
  \,.
\end{align*}
By comparing \eqref{eq:Pe+}, \eqref{eq:Pe-} with \eqref{eq:D0nija}, \eqref{eq:Dnij}, we obtain
\begin{align}
  \label{eq:P=}
  &\dP_n^+ = \widetilde \dD_n^0
  &&\dP_n^- = \dD_n
  \,.
\end{align}
The second equality in \eqref{eq:p+-} follows from \eqref{eq:Dn}, \eqref{eq:D0n}.
\qed

\medno
\section{Proof of Theorem \ref{thm:main}}
\label{sec:proof}

\noindent
In this section we prove Theorem \ref{thm:main}
following the strategy described in \cite[Sect. 3]{cesiOctopus}.

Since $G$ is a finite group, we can always assume
the representations are unitary with respect to some 
(positive definite) inner product $\<\cdot, \cdot\>$ defined
on the representation space $V$.
This will ensure that, 
if $w$ is a symmetric element of the group algebra, then
$\D_G(w,\dR)$ is self-adjoint.
For a self-adjoint linear operator $A$ we write $A\ge 0$ if 
$\<A \cdot, \cdot\>$ is
a positive semidefinite bilinear form.
We will write $\D_G(w,\dR) \ge 0$ if $\D_G(w,\dR)$ is positive
semidefinite for some (equivalently for each) unitary version
of $\dR$. 
We can thus define
\begin{align}
  \G(G) &= 
  \{ w \in \bC G^{(s)} : 
  \D_{G}(w, \dR) \ge 0, 
  \ \forall \dR\in \rep(G)  
  \}
  \,.
\end{align}
For future reference we summarize a few elementary properties
of the set $\G(G)$ in the following proposition.

\begin{proposition}\label{thm:gam}
We have:
\begin{enumerate}[(1)]
\addtolength\itemsep{1mm}
\item 
$\G(G)$ is a convex cone, \ie if $w,z \in \G(G)$, then
for any $\a,\b\in \bR_+$, $\a w+\b z \in \G(G)$;
\item 
$\D_G(w, \oplus_{i=1}^n \dR_i ) \ge0$ if and only if 
$\D_G(w,\dR_i) \ge0$ for every $i=1,\ldots,n$;
\item
  $\G(G)=
  \{ w \in \bC G^{(s)} : 
  \D_{G}(w, \dT) \ge 0, 
  \ \forall \dT\in \Irr(G)  
  \}
  $;
\item
  $\G(G)=
  \{ w \in \bC G^{(s)} : 
  \D_{G}(w, \dL ) \ge 0
  \}
  $, where $\dL$ is the left regular representation of $G$;
\item
  $  \bR_+ G^{(s)} \sset \G(G)$;
\item
  If $H$ is a subgroup of $G$, then $\G(H) \sset \G(G)$.
\end{enumerate}
\end{proposition}

\Pro\ 
(1) and (2) follow from the definitions.
(3) follows from (2) and \eqref{eq:Y}.
(4) follows from (2), (3) and \eqref{eq:LT}.
If $\dR$ is a unitary
representation on $V$, and $w\in \bC G^{(s)}$,
a straightforward computation
(see Proposition 2.1 in \cite{cesiOctopus})
yields 
\begin{align*}
  & \bigl\< \D_G(w,\dR) v, v \bigr\> =
  \ov2 \sum_{g\in G} w_g \, \|\dR(g) v - v\|^2
  &&v\in V
  \,.
\end{align*}
Thus $\D_G(w,\dR) \ge 0$ if $w$ is positive, which proves (5).

Finally, let $w\in \G(H)$. Then $\D_H(w,\dT) \ge 0$ for every
$\dT \in \Irr(H)$. If $\dS\in \Irr(G)$, then we have a
branching rule
\begin{align*}
  \dS \resR^{G}_{H} \iso \bigoplus_{\dT \in\Irr(H)} k(\dT) \, \dT
\end{align*}
where $k(\dT)$ are suitable nonnegative integers.
Since $w$ is an element of the group algebra of $H$,
the same decomposition applies to the representation Laplacian
\begin{align*}
  \D_{G}( w, \dS) = 
  \bigoplus_{\dT \in\Irr(H)} k(\dT) \, \D_{H}(w, \dT)
  \,.
\end{align*}
Therefore $\D_G(w,\dS)\ge 0$ and (6) follows.
\qed

\medno
In the following we regard $W_{n-1}$ as the subgroup of $W_n$
which fixes the last coordinate, that is
\begin{align*}
  W_{n-1} \iso \{ (\h,\pi) \in W_n : \h_n=0 \text{ and }
  \pi(n) = n
  \}
  \,.
\end{align*}
The key point of the proof is the following ``semirecursive''
result:

\begin{proposition}\label{thm:semi}
Let $w\in \bR_+ W_n^{(s)}$ and $z\in \bR_+ W_{n-1}^{(s)}$, be such that
$w-z \in \G(W_n)$. Then
\begin{equation}
  \psi_{W_n}(w) \ge \min\bigl\{
    \psi_{W_{n-1}}(z), \, 
    \psi_{W_n}(w, \dP_n)
  \bigr\}
  \,.
\end{equation}
\end{proposition}

\Pro\ 
Let $\cI_n$ be the set of all irreducible representations of $W_n$
that, when restricted to $W_{n-1}$, contain the trivial
representation.
The branching rule \eqref{eq:bra}
implies that
\begin{align}
  \cI_n
  =
  \bigl\{
  \dI_n,
  \ 
  [(n-1,1), \emp],
  \ 
  [(n-1),(1)]
  \bigr\}
  \,.
\end{align}
Thanks to Proposition \ref{thm:Pn} and \eqref{eq:Y2}, and
using the fact that $\psi_G(w,\dI) = +\oo$, we obtain
\begin{align*}
  \psi_{W_n}(w, \dP_n) = \min\{ 
    \psi_{W_n}(w, [(n-1,1),\emp] ),
    \psi_{W_n}(w, [(n-1), (1)] )
  \}
  \,.
\end{align*}
Thus Proposition \ref{thm:semi} follows from Proposition 3.2
in \cite{cesiOctopus}.
\qed

\medno
Let $\cA_n$ be the subset of $\bR_+ W_n^{(s)}$ considered
in the hypothesis of Theorem \ref{thm:main}
\begin{align}
  \label{eq:An}
  \cA_n = 
  \Bigl\{ 
  w = \sum_{i=1}^n a_i\, s_{\{i\}} + \sum_{(ij)\in T_n} b_{ij}\, (ij)
  : a_i\ge0,\; b_{ij} \ge 0
  \Bigr\}
  \,.
\end{align}
If $w\in \cA_n$, let us write $w = w_N + w_T$ with
\begin{align}
  \label{eq:wnt}
  &w_N = \sum_{i=1}^n a_i\, s_{\{i\}} 
  &&w_T = \sum_{(ij)\in T_n} b_{ij}\, (ij)
  \,.
\end{align}
We observe that in the $(\h,\pi)$ notation for the elements
of $W_n$ we have 
\begin{align}
  \label{eq:si}
  s_{\{i\}} = \bigl( \h_{\{i\}}, \idu_{\fS_n} \bigr)
  \quad\text{where}\quad
  (\h_{\{i\}})_j = \d_{ij} =
  \begin{cases}
  1 &\text{if $j = i$}\\
  0 &\text{if $j \ne i$.}
  \end{cases}
\end{align}
It follows from \eqref{eq:Dnij} and \eqref{eq:D0nija} that
\begin{align}
  \label{eq:D0si}
  \widetilde \dD_n^0\left(s_{\{i\}}\right) &= \id_n
  \\
  \label{eq:Dsi}
  \dD_n\left(s_{\{i\}}\right) &= \diag\left( (-1)^{\d_{ij}} \right)_{j=1}^n 
  \,.
\end{align}
Thus we get
\begin{align}
  \label{eq:w1}
  \D_{W_n}(w_N, \dont) &= 0
  \\
  \label{eq:w2}
  \D_{W_n}(w_N, \dD_n) &= 2 \diag( a_i)_{i=1}^n
  \\
  \label{eq:w3}
  \D_{W_n}(w_T, \dont) &= \D_{W_n}(w_T, \dD_n)
  \,.
\end{align}

\medno
\textit{Strategy for proving Theorem \ref{thm:main}}.

\smallno
Let us now assume that we find a map $\th : \cA_n \to \cA_{n-1}$ such that
the following holds for each $w\in \cA_n$:
\begin{enumerate}[(a)]
\item $w - \th(w) \in \G(W_n)$;
\item $\psi_{W_n}(w,\dP_n) \le \psi_{W_{n-1}}(\th(w), \dP_{n-1})$. 
\end{enumerate}
Then we can prove \eqref{eq:main} by induction.
Assume in fact that \eqref{eq:main} holds for $n=k-1$, that is
\begin{align}
    \label{eq:ind}
    &\psi_{W_{k-1}}(z) = \psi_{W_{k-1}}(z, \dP_{k-1} )
    &&\forall z\in \cA_{k-1}
    \,.
\end{align}
From Proposition \ref{thm:semi} and \eqref{eq:ind} with $z=\th(w)$,
and from properties (a), (b) of the map $\th$ it follows that
\begin{align*}
  \psi_{W_k}(w) \ge \min\bigl\{
    \psi_{W_{k-1}}( \th(w), \dP_{k-1}), 
    \psi_{W_{k}}( w, \dP_{k})
  \bigr\} 
  = \psi_{W_k}(w,\dP_k)
  \,,
\end{align*}
which, combined with the reversed inequality which is a trivial
consequence of \eqref{eq:minrep1}, implies 
$\psi_{W_k}(w) = \psi_{W_k}(w,\dP_k)$.
The induction step is completed.

In the next proposition we take care of the
starting point of the induction, $n=2$.

\begin{proposition}\label{thm:w2}
If $w\in \cA_2$, then $\psi_{W_2}(w) = \psi_{W_2}(w,\dP_2)$.
\end{proposition}

\Pro\
We have
\begin{align}
  \label{eq:irrw2}
  \Irr(W_2) = \{ \dI_2, \; [ (1,1), \,\emp ],\;
  [(1),\,(1)],\;
  [ \emp,\, (2) ],\;
  [\emp,\,(1,1) ]
  \}
  \,.
\end{align}
Proposition II.1 of \cite{GeiKin} states that if $[\a,\b] \in \Irr(W_n)$,
then
\begin{align}
  \label{eq:ba}
  [\b,\a] \iso \dJ_n \,\itensor \,[\a,\b]
  \,,
\end{align}
where $\itensor$ denotes the inner tensor product of representations
and $\dJ_n = [\emp,(n)]$ is the one-dimensional representation of $W_n$
given by
\begin{align}
  \label{eq:Jn}
  &\dJ_n( \h, \pi ) = (-1)^{\card\{ i \in \{1,\ldots,n\} \,:\, \h_i = 1 \} }
  &&(\h,\pi)\in W_n
  \,.
\end{align}
Using \eqref{eq:emp} and \eqref{eq:si},
we have
\begin{align*}
  \dT^{(\emp,\a)}( s_{\{i\}} ) 
  &= \dJ_n( \h_{\{i\}},\idu_{\fS_n} ) \, \dT^\a( \idu_{\fS_n} ) 
  = (-1) \, \dT^\a( \idu_{\fS_n} ) 
  = - \id_d 
  \\
  \dT^{(\emp,\a)}( (ij) ) 
  &= \dJ_n( 0, (ij) ) \, \dT^\a( (ij) ) 
  = \dT^\a( (ij) )
  \,,  
\end{align*}
where $d$ is the degree of $[\a]$. This
implies that, for every $w\in \cA_n$ of the form \eqref{eq:An}, we have
\begin{align}
  \label{eq:0a}
  \D_{W_n}(w, [\emp,\a] ) -
  \D_{W_n}(w, [\a,\emp] ) 
  =
  \bigl( 2 \sum_{i=1}^n a_i \bigr) \, \id_d
  \,. 
\end{align}
Therefore the eigenvalues of $\D_{W_n}(w, [\emp,\a] )$
are shifted, with respect to the eigenvalues of $\D_{W_n}(w, [\a,\emp] )$
by a nonnegative quantity. 
In particular, if $\a\ne(n)$, then
$[\a,\emp]$ is nontrivial and it has a spectral gap which is not greater
than the spectral gap of $[\emp,\a]$.
For this reason, representations of type $[\emp, \a]$ 
with $\a\ne (n)$ can
be safely omitted in the minimization process \eqref{eq:minrep2}
which produces the spectral gap of $w$.

Going back to the case $n=2$, we can take care of 
the representation $[\emp,(2)]$ with an explicit calculation.
If $w\in \cA_2$, it can be written as
\begin{align*}
  &w = x \, s_{\{1\}} + y \, s_{\{2\}} + z \, (12)
  &&x,y,z \ge 0
  \,.
\end{align*}
Since the Laplacian of the trivial representation is null, 
\eqref{eq:0a} becomes
\begin{align*}
  \D_{W_2}(w, [\emp,(2)] ) = 2 (x+y) \, \id_1
  \,.
\end{align*}
On the other hand, using \eqref{eq:Dnij}, we get
\begin{align*}
  \D_{W_2}(w,\dD_2) =
  \begin{bmatrix}
  2x + z & -z \\ -z & 2y+z
  \end{bmatrix}
  =:B_2
\end{align*}
with spectral gap
\begin{align*}
  \psi_{W_2}(w,\dD_2) &= \l_1(B_2) = 
  x + y + z - ( (x-y)^2 + z^2 )^{1/2}
  \\
  &\le x+y \le 2(x+y) = \psi_{W_2}(w,[\emp,(2)])
  \,.
\end{align*}
Thus, for the purpose of computing the spectral
gap of $w$, representation $[\emp,(2)]$ can also be disregarded in
the list \eqref{eq:irrw2}. By consequence
\begin{align*}
  \psi_{W_2}(w) = \min\{ 
    \psi_{W_2}(w, [(1,1),\emp]), \;
    \psi_{W_2}(w, [(1),(1)])
  \}
  = \psi_{W_2}(w,\dP_2)
  \,.
\qed
\end{align*}

\medno
\textit{The mapping $\th$}.
In order to conclude the proof of Theorem \ref{thm:main}
we are going to define a map $\th :\cA_n \to \cA_{n-1}$
which satisfies properties (a) and (b) stated above.

\smallno
If $w = \sum_{i=1}^n a_i\, s_{\{i\}}$ with $a_i\ge 0 $, we let
$\ell$ be the largest index $j$ such that 
$a_j = \min_i a_i$,
and we define
\begin{align}
  \label{eq:thn}
  \th^N(w) = \sum_{ i=1,\; i\ne \ell }^n  a_i \, s_{\{i\}}
  \,.
\end{align}

\smallno
If $w = \sum_{(ij)\in T_n} b_{ij}\, (ij)$ with $b_{ij}\ge 0$,
for each $m=1,\ldots, n$ we let
\begin{align}
  \label{eq:capu}
  \th^T_m (w) = 
  \sum_{ \atopp{(ik) \in T_{n}}{  i,k\ne m }} 
  \left[ w_{ik} + \frac{w_{im} \, w_{k m}}
  {\sum_{j\ne m} w_{j m}}
  \right] \, (ik)
  \,.
\end{align}
Finally we define a mapping $\th(w)$ as follows:
let $w = w_N + w_T$ as in \eqref{eq:wnt}.
Then we let
\begin{align}
  \label{eq:th}
  \th(w) = \th^N(w_N) + \th^T_\ell(w_T)
\end{align}
where $\ell$ is defined as above.

\begin{remark}\label{thm:ass}
We can assume,
without loss of generality, that $\ell = n$
in \eqref{eq:th}.
In this way $\th (\cA_n) \sset \cA_{n-1}$.
Otherwise one can define $W_{n-1}$ as the subgroup of $W_n$
obtained by ``dropping
the $\ell^{th}$ coordinate''.
\end{remark}

\begin{remark}
The mapping $\th^T_m$, amazingly, 
appeared almost simultaneously
in the preprint versions of \cite{Diek} and \cite{CaLiRi}
and it was key point which, together with a quite tricky
inequality, called the ``octopus inequality'', produced
a proof of Aldous' spectral gap conjecture in \cite{CaLiRi}.
\end{remark}

\medno
Properties (a) and (b) of the mapping $\th$ will be
proved in Lemmas \ref{thm:octo} and \ref{thm:P} respectively,
completing in this way the proof of Theorem \ref{thm:main}.

\begin{lemma}\label{thm:octo}
If $w \in \cA_n$, then $w - \th(w) \in \G(W_n)$.
\end{lemma}

\Pro\ 
We can write
\begin{align*}
  w - \th(w) = \underbrace{w_N - \th^N(w_N)}_{\d w_N} 
  + \underbrace{w_T - \th^T(w_T)}_{\d w_T}
  \,,
\end{align*}
Since $\d w_N$ is positive and symmetric, we have $\d w_N \in \G(W_n)$,
thanks to Proposition \ref{thm:gam}(5).

On the other hand Theorem 2.3 of \cite{CaLiRi},  
the ``octopus inequality'' (see also Section 4 of \cite{cesiOctopus}
for a slightly simpler proof in which the algebraic
perspective is more explicit),
states that
$\d w_T \in \G(\fS_n)$.
By Proposition \ref{thm:gam}(6) we get
$\d w_T \in \G(W_n)$.

Hence $w-\th(w) = \d w_N + \d w_T \in \G(W_n)$
by Proposition \ref{thm:gam}(1)
\qed

\medno

\begin{lemma}\label{thm:P}
If $w \in \cA_n$, then
\begin{align}
  \psi_{W_n}( w , \dP_n ) \le \psi_{W_{n-1}}( \th w, \dP_{n-1} )
  \,.
\end{align}
\end{lemma}

\Pro\ 
From Proposition \ref{thm:Pn} we know that $\dP_n^- = \dD_n$
is irreducible, while $\dP_n^+ = \widetilde \dD_n^0$ contains
the trivial representation with multiplicity $1$, hence, by \eqref{eq:gapt},
we obtain
\begin{equation}
\begin{split}
  \label{eq:gapw}
  \psi_{W_n}( w , \dP_n ) 
  &= \min \{ \psi_{W_n}( w , \dP^+_n ),
    \psi_{W_n}( w , \dP^-_n )
  \}
  \\
  &=
  \min\{ 
    \l_2\left( \D_{W_n}(w, \dP^+_n) \right),
    \l_1\left( \D_{W_n}(w, \dP^-_n) \right)
  \}
  \,.
\end{split}
\end{equation}
Since $\supp(\th w) \in W_{n-1}$, the last row and column
of its Laplacian are zero, thus we can write its representation
Laplacian in block diagonal form as
\begin{align}
  \label{eq:n-1}
  \D_{W_n}( \th w, \dP^\pm_n ) =
  \D_{W_{n-1}}( \th w, \dP^\pm_{n-1} ) \oplus [0]_{1\times 1}
  \,,
\end{align}
where $[x]_{1\times 1}$ is the $1\times 1$ matrix whose
unique entry is equal to $x$.
This implies 
\begin{equation}
\begin{split}
  \label{eq:gapwth}
  &\psi_{W_{n-1}}( \th w , \dP_{n-1} ) 
  \\
  &\quad=
  \min\{ 
    \l_2\left( \D_{W_{n-1}}(\th w, \dP^+_{n-1}) \right),
    \l_1\left( \D_{W_{n-1}}(\th w, \dP^-_{n-1}) \right)
  \}
  \\
  &\quad=
  \min\{ 
    \l_3\left( \D_{W_{n}}(\th w, \dP^+_{n}) \right),
    \l_2\left( \D_{W_{n}}(\th w, \dP^-_{n}) \right)
  \}
  \,.
\end{split}
\end{equation}
We write $w = w_N + w_T$ with $w_N$ and $w_T$ as in \eqref{eq:wnt}.
For simplicity we also define the following matrices:
\begin{align*}
  M_n &= \D_{W_n}(w_T, \dP^+_n)
  &
  M^\th_n &= \D_{W_n}( \th^T  w_T, \dP^+_n)
  \\
  F_n &= 2 \diag( a_i)_{i=1}^n
  &
  F_n^\th &= 2 \diag( a_1, \ldots, a_{n-1}, 0 )
  \,.
\end{align*}
We are assuming (remember Remark \ref{thm:ass})
the $a_n = \min_j a_j$.
It follows from \eqref{eq:P=}, \eqref{eq:w1}, \eqref{eq:w2}, \eqref{eq:w3}
that
\begin{align}
  \label{eq:M}
  \D_{W_n}(w, \dP^+_n) &= M_n
  &
  \D_{W_n}(\th w, \dP^+_n) &= M^\th_n
  \\
  \label{eq:MF}
  \D_{W_n}(w, \dP^-_n) &= M_n + F_n
  &
  \D_{W_n}(\th w, \dP^-_n) &= M^\th_n + F^\th_n
  \,.
\end{align}
By \eqref{eq:n-1} and \eqref{eq:MF} we can write
\begin{equation}
  \label{eq:n-1a}
  M^\th_n + F^\th_n = B_{n-1} \oplus [0]_{1\times 1}
\end{equation}
with $B_{n-1} = \D_{W_{n-1}}( \th w, \dP^-_{n-1} )$.
But then we have 
\begin{align*}
  M^\th_n + F_n = B_{n-1} \oplus [2 a_n]_{1\times 1}
  \,.
\end{align*}
By consequence%
\footnote{the spectrum is always considered as a multiset, so if,
for instance, $\spec(A) = \{ 0, 0, 1, 2 \}$, then 
$(\spec(A)\setm \{0\})\cup\{1\} = \{ 0, 1, 1, 2 \}$.}
\begin{align}
  \label{eq:specMF}
  \spec( M^\th_n + F_n ) = 
  \spec( B_{n-1} ) \cup \{2 a_n\}
  \,.
\end{align}
Since $\th^T w_T$ is symmetric and positive, by 
Proposition \ref{thm:gam}(5) 
the matrix $M^\th_n$
is positive semidefinite, which implies (see, for instance,
\cite[Corollary 4.3.3]{HorJohn})
\begin{align}
  \label{eq:lk}
  &\l_k(M^\th+ F_n) \ge \l_k(F_n)
  &&k=1,\ldots, n
  \,.
\end{align}
Thus we get
\begin{align*}
  2 a_n = \l_1( F_n) \le \l_1( M^\th_n + F_n)
  \,.
\end{align*}
But \eqref{eq:specMF} says that $2 a_n$ is actually an eigenvalue
of $M^\th_n+F_n$, so it must be the lowest one
\begin{align}
  \label{eq:2a}
  2 a_n = \l_1( M^\th_n + F_n)
  \,.
\end{align}
Therefore, by \eqref{eq:n-1a}, we get
\begin{align}
  \label{eq:l22}
  \l_2(M^\th_n + F_n) = \l_1(B_{n-1}) = \l_2(M^\th_n + F^\th_n)
  \,.
\end{align}
Using the explicit expression \eqref{eq:D0nija} for the matrix elements
of the representation $\dP^+_n = \widetilde \dD^0_n$,
it is straightforward to check that the matrix elements of 
$M_n - M_n^\th$ are given by
\begin{equation}
  \label{eq:dM}
  [L_n]_{ij} := \bigl[ M_n - M_n^\th \bigr]_{ij} = \frac{d_i \,d_j}{d_n} \,,
\end{equation}
where $d_i = -b_{in}$ for $i=1,\ldots,n-1$ and $d_n = \sum_{i=1}^{n-1} b_{in}$.
Following \cite{Diek} we observe
that $L_n$
is a rank-1 matrix, so by standard linear algebra results
as \cite[Thm. 4.3.4]{HorJohn}, 
one obtains, in particular, that 
\begin{align}
  \label{eq:l22a}
  \l_2(M_n ) &= \l_2( M^\th_n + L_n) \le \l_3(M^\th_n)
  \\
  \label{eq:l12}
  \l_1(M_n+F_n) &= \l_1( M^\th_n+F_n + L_n ) \le \l_2(M^\th_n+F_n)
  \,.
\end{align}
Thus, using \eqref{eq:l22}, we get
\begin{align}
  \label{eq:lth}
  \l_1(M_n+F_n) \le \l_2( M^\th_n + F^\th_n)
  \,.
\end{align}
Relations \eqref{eq:gapw}, \eqref{eq:gapwth}, \eqref{eq:M}, \eqref{eq:MF},
\eqref{eq:l22a} and \eqref{eq:lth} imply Lemma \ref{thm:P}
\qed

\medno
\section{A few concluding remarks and one open problem}
\label{sec:rem}

\noindent
Theorem \ref{thm:main}, together with Proposition \ref{thm:Pn}, 
states that if $w\in \cA_n$, then
the representation ``responsible'' for the spectral gap
is either $\dD_n \iso [(n-1),(1)]$ or $\widetilde \dD_n^0$
that is $[(n-1,1),\emp]$, since the trivial summand in \eqref{eq:p+-}
plays no role. 
In our proof we are led to consider these two representations
because they are the representations which, when restricted to $W_{n-1}$
contain the trivial one.

We show that this is not an artifact of our strategy:
we actually need to include both of them, that is
the statement of Theorem \ref{thm:main}
cannot be strengthened by replacing $\dP_n$ with
either $\dD_n$ or $\widetilde\dD_n^0$.
Let 
\begin{align*}
  w_N &= \sum_{i=1}^n a_i\, s_{\{i\}}  
  &
  w_T &= \sum_{(ij)\in T_n} b_{ij}\, (ij)
  && a_i \ge0, \ b_{ij}\ge 0
  \,.
\end{align*}
Let $\uu a := \min_i  a_i$ and assume $\uu a > 0$.
Assume also that:
\begin{enumerate}[(i)]
\item 
there are enough strictly positive $b_{ij}$
so that $\supp(w)$ generates $S_n$.
\end{enumerate}
This condition is equivalent to requiring that the graph on
$\{1,\ldots,n\}$ with edge set
$\cE = \{  \{i,j\} : b_{ij}>0 \}$ is connected.

\smallno
For $\e>0$ define
\begin{align*}
  w_\e = w_N + \e \, w_T
  \,.
\end{align*}
Thanks to \eqref{eq:w1}, \eqref{eq:w2}, \eqref{eq:w3}
we can write
\begin{align*}
  F_n &:= \D_{W_n}(w_N, \dD_n) = 2 \diag( a_i)_{i=1}^n
  \\
  M_n &:= \D_{W_n}(w_T, \dD_n) = \D_{W_n}(w_T, \dont)
  \,.
\end{align*}
Hence
\begin{align*}
  \D_{W_n}(w_\e, \dD_n) &= 
  F_n + \e\, M_n
  \\
  \D_{W_n}(w_\e, \dont) &= 
  \e\, M_n
  \,.
\end{align*}
The lowest eigenvalue of $M_n$ is trivially $\l_1(M_n)=0$ with eigenvector
equal to any constant vector.
It is easy to show that hypothesis (i) above
implies that $0$ is a simple eigenvalue, that is $\l_2(M_n) > 0$
(see, for instance, Proposition 2.1 of \cite{cesiOctopus}).
By perturbation theory we obtain, using \eqref{eq:gapt},
\begin{align*}
  \psi_{W_n}(w_\e, \dD_n) &= \l_1( F_n + \e M_n ) =  2 \uu a + \cO(\e)
  \\
  \psi_{W_n}(w_\e, \dont) &= \e \l_2( M_n ) =  \cO(\e)
  \,,
\end{align*}
where $\cO(\e)$ is a generic quantity which goes to $0$ as $\e\to 0^+$.
Hence, for small $\e$ the spectral gap of $w_\e$
is determined by $\dont$.

\medno
Consider now the opposite situation with
\begin{align*}
  w_\e = \e \, w_N + \, w_T
  \,.
\end{align*}
We obtain
\begin{align*}
  \psi_{W_n}(w_\e, \dD_n) &= \l_1( M_n ) + \cO(\e) = \cO(\e)
  \\
  \psi_{W_n}(w_\e, \dont) &= \l_2( M_n ) > 0
  \,,
\end{align*}
hence, in this case, for $\e$ small enough, the spectral gap of $w_\e$
is determined by $\dD_n$.

\medno
Lastly we want to discuss the possibility
of proving our main theorem 
for more general elements $w$ of the group algebra
than those considered in \eqref{eq:wmain}.
For $A\sset \{1,\ldots,n\}$,
let $s_A$ be the element of $W_n$ which
in the defining representation is given by the diagonal matrix
\begin{align}
  \dD_n(s_A) = \diag( x_i )_{i=1}^n
  \quad
  \text{where}
  \quad
  x_i = \begin{cases}
  -1 & \text{if $i\in A$} \\
  +1 & \text{if $i\notin A$.} 
  \end{cases}
\end{align}
Let then $Y^+_n$ ($Y^-_n$) be the set of all subsets of $\{1,\ldots,n\}$
of even (odd) cardinality, and let
\begin{align*}
  w_N^\pm &= \sum_{A\in Y^\pm_n} a_A s_A 
  &
  w_T &= \sum_{(ij)\in T_n} b_{ij}\, (ij)
  && a_A \ge0, \ b_{ij}\ge 0
  \,.
\end{align*}

\begin{question}
\label{th:que}
Does the equality $\psi_{W_n}(w) = \psi_{W_n}(w, \dP_n )$
also hold for elements $w$ of the form $w = w^+_N + w^-_N + w_T$?
\end{question}

\noindent
We show that the answer is (in general) negative.

\smallno
Let
\begin{align*}
  \hat a^\pm_i = \sum_{A \in Y^\pm_n \,:\, A \ni i} a_A 
\end{align*}
and assume that:
\begin{enumerate}[(i)]
\addtolength\itemsep{0.5mm}
\item 
there are enough strictly positive $b_{ij}$
so that $\supp(w)$ generates $S_n$;
\item
$w_N^- \ne 0$;
\item 
$\uu a^+ := \min_i \hat a^+_i > 0$.
\end{enumerate}
The first two conditions are necessary, since otherwise
the support of $w$ does not generate $W_n$,
the spectral is trivially zero, and the problem becomes uninteresting.

\smallno
For $\e>0$, let 
\begin{align}
  \label{eq:we}
  w_\e = w_{N}^+ + \e \, w_{N}^- + w_T
  \,.
\end{align}
Thanks to \eqref{eq:ba} and \eqref{eq:Jn} we have
\begin{align*}
  &\dT^{\b,\a}( w_{N}^+ ) 
  = 
  \dT^{\a,\b}( w_{N}^+ ) 
  &&
  \dT^{\b,\a}( w_T ) 
  = 
  \dT^{\a,\b}( w_T ) 
\end{align*}
which, since $[(n),\emp]$ is the trivial representation $\dI_n$, 
implies in particular that
\begin{align}
  \label{eq:D0}
  \D_{W_n}( w_{N}^+ , [\emp,(n)] )  
  = 
  \D_{W_n}( w_T, [\emp,(n)] )  
  = 
  0
  \,.
\end{align}
As for the ``odd term'' $w_N^-$, using \eqref{eq:Jn} we obtain
\begin{align*}
  \dT^{(\emp,\a)}( s_A ) 
  = (-1) \, \dT^\a( \idu_{\fS_n} ) 
  = - \id_d 
  \,,  
\end{align*}
where $d$ is the degree of the representation $[\a]$, thus
\begin{align}
  \label{eq:De}
  \D_{W_n} (w_N^-, [\emp,\a] ) =
  2 \Bigl(\sum\nolimits_{A\in Y_n^-} a_A \Bigr) \id_d
  \,.
\end{align}
From \eqref{eq:D0} and \eqref{eq:De} it follows that
\begin{align*}
  \psi_{W_n}( w_\e, [\emp,(n)] ) 
  =
  \e \l_1\bigl( \D_{W_n}( w_N^-, [\emp,(n)] ) \bigr)
  = 2 \e \sum\nolimits_{A\in Y_n^-} a_A
  = \cO(\e)
  \,.
\end{align*}
On the other hand we claim that $\psi_{W_n}(w_\e, \dP_n)$
can be bounded from below by a strictly positive (independent of $\e$)
quantity. 
It easy to see that equalities \eqref{eq:w1} and \eqref{eq:w2} become
\begin{align}
  \label{eq:w1a}
  \D_{W_n}(w_N^\pm, \dont) &= 0
  \\
  \label{eq:w2a}
  \D_{W_n}(w_N^\pm, \dD_n) &= 2 \diag( \hat a^\pm_i)_{i=1}^n
  \,.
\end{align}
Let (remember \eqref{eq:w3})
\begin{align*}
  F_n^\pm &:= \D_{W_n}(w_N^\pm, \dD_n) = 2 \diag( \hat a^\pm_i)_{i=1}^n
  \\
  M_n &:= \D_{W_n}(w_T, \dD_n) = \D_{W_n}(w_T, \dont)
  \,.
\end{align*}
In this way we have obtained
\begin{align*}
  \D_{W_n}( w_\e, \dD_n ) &= M_n + F^+_n + \e F^-_n \ge F^+_n 
  \\
  \D_{W_n}( w_\e, \dont ) &= M_n
  \,,
\end{align*}
where the inequality is intended in the sense of quadratic forms.
Thanks to assumption (i) above, we know that $\l_2(M_n)$
is strictly positive, 
therefore
\begin{align*}
  \psi_{W_n}( w_\e, \dP_n) 
  &=
  \min\{ \psi_{W_n}( w_\e, \dD_n), \psi_{W_n}( w_\e, \dont) \}
  \\
  &\ge \min\{ \l_1(F_n^+) , \l_2(M_n) \} 
  =
  \min\{ 2 \uu a^+ , \l_2(M_n) \} > 0
  \,.
\end{align*}
Thus, for $\e$ small enough, we have 
$\psi_{W_n}( w_\e, [\emp,(n)] ) < \psi_{W_n}( w_\e, \dP_n)$,
which implies a negative answer to Question \ref{th:que}.

We observe that a crucial element for this ``counterexample''
is assumption (iii) above. This leaves room for
a conjecture.%

\begin{conjecture}
\label{th:wro}
If $w = w_N^- + w_T$, then $\psi_{W_n}(w) = \psi_{W_n}(w, \dP_n)$.
\end{conjecture}

\noindent
The most obvious approach for proving this result would be
to generalize the map $\th^N$ of \eqref{eq:thn} as
\begin{align}
  \label{eq:thn0}
  \th^N(w_N^-) = \sum_{ A\in Y^-_n\,:\, A\not\ni \ell}  a_A s_A
  \,,
\end{align}
where $\ell$ is the largest index $j$ such that
$\hat a_j^- = \min_i \hat a_i^-$,
Unfortunately this does not work because, with this choice,
Lemma \ref{thm:P} is false. A counterexample 
can be found already for $n=3$: if
\begin{equation*}
  w = \sum_{i=1}^3 s_{\{i\}} + s_{\{1,2,3\}} + (12) + (23) + (13) 
\end{equation*}
then
\begin{equation*}
  \th(w) = s_{\{1\}} + s_{\{2\}} + \frac{3}{2} \, (12) 
\end{equation*}
which produces $\psi_{W_3}(w, \dP_3) = 3 > 2 = \psi_{W_2}(\th(w), \dP_2)$.
We emphasize that this is a counterexample to Lemma \ref{thm:P},
\textit{not} to Conjecture \ref{th:wro}, since we have in fact
$\psi_{W_3}(w) = 3$.
Hence one should devise a different map $\th : \cA_n \to \cA_{n-1}$,
keeping in mind that there is a delicate balance between
the two properties (a) and (b) of Section \ref{sec:proof}
which must be satisfied by $\th$.

\medno
\textbf{Acknowledgements.}
In the first version of this paper we erroneously claimed
to have proven Conjecture \ref{th:wro} due to a 
a mistake in the proof of (a more general version of) 
Lemma \ref{thm:P} where \eqref{eq:thn0} was used.
We thank one of the referees for finding the mistake
in the proof, which prompted us to find the above
counterexample.


\noindent
{
\hbadness=10000

\begin{thebibliography}{FOW85}

\bibitem[Ald]{Ald}
D.~Aldous,
  \emph{www.stat.berkeley.edu/\%7{E}aldous/{R}esearch/{O}{P}/sgap.html}.

\bibitem[Bac94]{Bac}
R.~Bacher, \emph{Valeur propre minimale du laplacien de {C}oxeter pour le
  groupe sym\'etrique}, J. Algebra \textbf{167} (1994), no.~2, 460--472.

\bibitem[Ces09]{Ces3}
F.~Cesi, \emph{Cayley graphs on the symmetric group generated by initial
  reversals have unit spectral gap}, Electron. J. Combin. \textbf{16} (2009),
  no.~1, Note 29, 7.

\bibitem[Ces10]{Ces1}
\bysame, \emph{On the eigenvalues of {C}ayley graphs on the symmetric group
  generated by a complete multipartite set of transpositions}, J. Algebraic
  Combin. \textbf{32} (2010), no.~2, 155--185.

\bibitem[Ces16]{cesiOctopus}
\bysame, \emph{A few remarks on the octopus inequality and {A}ldous' spectral
  gap conjecture}, Communications in Algebra \textbf{44} (2016), no.~1,
  279--302.

\bibitem[CLR10]{CaLiRi}
P.~Caputo, T.~M. Liggett, and T.~Richthammer, \emph{Proof of {A}ldous' spectral
  gap conjecture}, J. Amer. Math. Soc. \textbf{23} (2010), no.~3, 831--851.

\bibitem[CT17]{chung2017spectral}
F.~Chung and J.~Tobin, \emph{The spectral gap of graphs arising from substring
  reversals}, The Electronic Journal of Combinatorics \textbf{24} (2017),
  no.~3, 3--4.

\bibitem[Die10]{Diek}
A.~B. Dieker, \emph{Interlacings for random walks on weighted graphs and the
  interchange process}, SIAM J. Discrete Math. \textbf{24} (2010), no.~1,
  191--206.

\bibitem[DS81]{DiSh}
P.~Diaconis and M.~Shahshahani, \emph{Generating a random permutation with
  random transpositions}, Z. Wahrsch. Verw. Gebiete \textbf{57} (1981), no.~2,
  159--179.

\bibitem[FOW85]{FOW}
L.~Flatto, A.~M. Odlyzko, and D.~B. Wales, \emph{Random shuffles and group
  representations}, Ann. Probab. \textbf{13} (1985), no.~1, 154--178.

\bibitem[GK78]{GeiKin}
L.~Geissinger and D.~Kinch, \emph{Representations of the hyperoctahedral
  groups}, Journal of algebra \textbf{53} (1978), no.~1, 1--20.

\bibitem[GP00]{GecPfe}
M.~Geck and G.~Pfeiffer, \emph{Characters of finite {C}oxeter groups and
  {I}wahori-{H}ecke algebras}, London Mathematical Society Monographs. New
  Series, vol.~21, The Clarendon Press Oxford University Press, New York, 2000.

\bibitem[HH17]{huang2017adjacency}
X.~Huang and Q.~Huang, \emph{The adjacency spectral gap of some cayley graphs
  on alternating groups}, \upshape arXiv:1711.08944 (2017).

\bibitem[HJ90]{HorJohn}
R.~A. Horn and C.~R. Johnson, \emph{Matrix analysis}, Cambridge University
  Press, Cambridge, 1990, Corrected reprint of the 1985 original.

\bibitem[HJ96]{HaJu}
S.~Handjani and D.~Jungreis, \emph{Rate of convergence for shuffling cards by
  transpositions}, J. Theoret. Probab. \textbf{9} (1996), no.~4, 983--993.

\bibitem[JK81]{JaKe}
G.~James and A.~Kerber, \emph{The representation theory of the symmetric
  group}, Encyclopedia of Mathematics and its Applications, vol.~16,
  Addison-Wesley Publishing Co., Reading, Mass., 1981.

\bibitem[Kas11]{Kas}
M.~Kassabov, \emph{Subspace arrangements and property {T}}, Groups, Geometry,
  and Dynamics \textbf{5} (2011), no.~2, 445--477.

\bibitem[KN97]{KoNa}
T.~Koma and B.~Nachtergaele, \emph{The spectral gap of the ferromagnetic {XXZ}
  chain}, Lett. Math. Phys. \textbf{40} (1997), no.~1, 1--16.

\bibitem[Mor08]{Mor}
B.~Morris, \emph{Spectral gap for the interchange process in a box}, Electron.
  Commun. Probab. \textbf{13} (2008), 311--318.

\bibitem[PP18]{parzanchevski2018aldous}
O.~Parzanchevski and D.~Puder, \emph{Aldous' spectral gap conjecture for normal
  sets}, \upshape arXiv:1804.02776 (2018).

\bibitem[SC11]{StCo}
S.~Starr and M.~P. Conomos, \emph{Asymptotics of the spectral gap for the
  interchange process on large hypercubes}, Journal of Statistical Mechanics:
  Theory and Experiment \textbf{2011} (2011), no.~10, P10018.

\end{thebibliography}

\def\polhk#1{\setbox0=\hbox{#1}{\ooalign{\hidewidth
  \lower1.5ex\hbox{`}\hidewidth\crcr\unhbox0}}}
\providecommand{\bysame}{\leavevmode\hbox to3em{\hrulefill}\thinspace}
\providecommand{\MR}{\relax\ifhmode\unskip\space\fi MR }
\providecommand{\MRhref}[2]{%
  \href{http://www.ams.org/mathscinet-getitem?mr=#1}{#2}
}
\providecommand{\href}[2]{#2}

}

\end{document}